
\input amstex
\documentstyle{amsppt}
\topmatter
\title Quantum Cohomology of Flag Varieties
\endtitle
\author Ionu\c t Ciocan-Fontanine
\endauthor
\affil University of Utah \endaffil

\address{Department of Mathematics, University of Utah, Salt Lake City,
UT 84112} \endaddress
\address{\rm email: ciocan-f\@math.utah.edu} \endaddress

\date{ May 2, 1995 } \enddate

\endtopmatter
\document
\heading Introduction \endheading

The quantum cohomology ring of a K\"ahler manifold $X$ is a deformation of the
 usual cohomology ring which appears naturally in theoretical physics in the
study of the supersymmetric nonlinear sigma models with target $X$. In \cite
{W}, Witten introduces the quantum multiplication of cohomology classes on
$X$ as a certain deformation of the usual cup-product, obtained by adding to it
 the so-called instanton corrections (see also \cite{V}). These can be in
turn interpreted as
intersection numbers on a sequence of moduli spaces of (holomorphic) maps
$\Bbb P^1\rightarrow X$. To make this interpretation rigorous according to
mathematical standards, one encounters severe problems, mainly because these
moduli spaces are not compact and they may have the wrong dimension.  Recently,
 substantial efforts have been made to put the theory on firm mathematical
footing, and a proof for the existence of the quantum cohomology ring, using
methods of symplectic topology, has been given by Ruan and Tian \cite{RT}, for
 a large class of manifolds ( semi-positive symplectic manifolds).

Nevertheless, computing the quantum cohomology ring is typically a difficult
task (the method of proving existence relies on changing the complex structure
 of $X$ to a generic almost complex structure, hence it is not suited for
computations). Several examples have been worked out (Batyrev (\cite {Bat})
for toric varieties - see also \cite {MP} - , Bertram (\cite {Be2}) and Siebert
-Tian (\cite {ST}) for Grassmannians), but the problem is far from solved.
In \cite {GK}, based on the conjectures following from conformal field theory,
 Givental and Kim proposed a presentation of the quantum cohomology ring in the
  case of flag varieties.

 In the present paper we describe a method for computing these rings, building
on the ideas in \cite {Be1}, \cite {Be2} and \cite {BDW}. We give a rigorous
construction of
 the (genus 0) Gromov-Witten invariants for the flag varieties (i.e. the
intersection numbers we mentioned above) and using it we complete our
computations, recovering the statement in \cite {GK}. Fulton's main result
in \cite {F1} will be essential for the proof. It is also worth noting that
the method  provides a new (algebraic-geometric!) proof for the existence of
quantum cohomology, as in \cite {Be2}. This will be done elsewhere. (The
proof in \cite {RT} has been recently redone in an algebraic setting for the
case of homogeneous spaces by J. Li and G. Tian - see \cite {LT}.)

Finally, we should remark that the method described here, coupled with the
results in \cite {F2}, should apply also to more general homogeneous spaces
(at least for the classical groups) giving a presentation for the quantum
cohomology ring in these cases as well.

{\bf Acknowledgement:} The results presented here will be part of my Ph.D.
thesis. I am extremely grateful to my advisor, Aaron Bertram, for introducing
me to the problem and for his timely and patient help.
\heading Notations and Statements of Results \endheading
We start by recalling some well-known facts about flag varieties and their
cohomology.

Let $V\cong \Bbb C^n$ be a complex vector space and define $F(n)=F(V^{\ast })$
to be the variety of complete flags:
 $U_{\bullet }:\{ 0\} =U_0\subset U_1\subset U_2\subset \dots \subset U_{n-1}
\subset U_n= V^{\ast }$.

 On $F(V^{\ast })$ there is an universal flag of subbundles $$E_1\subset E_2
\subset \dots \subset E_{n-1}\subset E_n=V^{\ast } \otimes \Cal O_{F(
V^{\ast })}$$ and an universal sequence of quotient bundles
$$ V^{\ast } \otimes
\Cal O_{F(V^{\ast })}=L_n\twoheadrightarrow L_{n-1}\twoheadrightarrow \dots
\twoheadrightarrow L_1 \twoheadrightarrow 0, $$ where
$L_i=V^{\ast } \otimes \Cal O_{F(V^{\ast })}/E_{n-i}$ for $i=1\dots n-1$.

Fix a complete flag of subspaces $V_{\bullet }^{\ast }: \{0 \} =V_0^{\ast }
\subset V_1^{\ast }\subset \dots \subset V_{n-1}^{\ast }\subset V_n^{\ast }=
V^{\ast }$. We have then induced maps $f_{pq}: V_p^{\ast } \otimes \Cal O_{F
(V^{\ast })}\rightarrow L_q$. Let $S_n$ be the symmetric group on $n$ letters.
For $w\in S_n$, let $r_w(q,p)=\text {card}\{ i\ |\ i\leq q, w(i)\leq p\} $.
Set $$\Omega _w=\Omega _w(V_{\bullet }^{\ast })=\{ U_{\bullet }\in F(V^{\ast })
\ |\ \text {rank}_{U_{\bullet }}f_{pq}\leq r_w(q,p), 1\leq p,q\leq n-1\} $$
and $$ X_w=X_w(V_{\bullet }^{\ast })=\{ U_{\bullet }\in F(V^{\ast })\ |\
\text {dim} (U_q\bigcap V_p^{\ast })\geq r_w(q,p), 1\leq p,q \leq n-1\} .$$
$\Omega _w$ is a subvariety of $F(V^{\ast })$ of codimension $\ell (w)=$
the length of the permutation $w\in S_n$. If we let $w_0\in S_n$ be
the permutation of longest length (given by $w_0(i)=n-i+1,\ i=1,\dots ,n$),
then $\Omega _w=X_{ww_0}$ for all $w\in S_n$

{\bf Facts: 1} $\{ X_w\} _{w\in S_n}$ and $\{\Omega _w\} _{w\in S_n}$ form
dual additive bases for $CH^{\ast }(F(V^{\ast }))\cong H^{\ast }(F(V^{\ast })
;\Bbb Z)$.

{\bf 2} Let $x_i=c_1(\text {ker}(L_i\rightarrow L_{i-1}))=c_1(
E_{n-i+1}/E_{n-i}),\ I=1,\dots ,n$. Then $\{ x_1^{i_1}x_2^{i_2}\dots
x_n^{i_n}\ |\ i_j\leq n-j \}$ form an additive basis as well, and
$$H^{\ast }(F(V^{\ast });\Bbb Z)\cong \Bbb Z[x_1,x_2,\dots x_n]/\ (R_1(n),
 R_2(n), \dots ,R_n(n)),$$
where $R_i(n)$ is the $i^{th}$ elementary symmetric function in
$x_1, x_2,\dots ,x_n$.

To express $\Omega _w$ in $H^{\ast }(F(V^{\ast });\Bbb Z)$ we need the notion
of Schubert polynomials (see \cite{LS1}, \cite{LS2} ). Define operators
$\partial _i,\ i=1,\dots ,n-1$ on $\Bbb Z[x_1,\dots ,x_n]$ by
$$\partial _iP=\frac {P(x_1,\dots ,x_n)-P(x_1,\dots ,x_{i-1},x_{i+1},x_i,
x_{i+2},\dots ,x_n)}{x_i-x_{i+1}}.$$
For any $w\in S_n$, write $w=w_0\cdot s_{i_1}\cdot \dots \cdot s_{i_k}$,
 with $k=\frac{n(n-1)}{2}-\ell (w)$, where $s_i=(i,i+1)$ is the transposition
interchanging $i$ and $i+1$. The polynomial $\Cal \sigma _w(x)\in \Bbb Z[x_1
,\dots ,x_n]$ defined by $$\Cal \sigma _w(x)=\partial _{i_k}\circ \dots
\circ \partial_{i_1}(x_1^{n-1}x_2^{n-2}\dots x_{n-1})$$ is the {\it Schubert
polynomial} associated to $w$. With this definition, we have a Giambelli type
formula:$$\Omega_w=\Cal \sigma_w(x_1,\dots ,x_n)$$ in $H^{\ast }(F(V^{\ast })
;\Bbb Z)$.

$F(n)$ embeds in a projective space by a Pl\"ucker embedding as a linear
section of a product of $n-1$ Grassmannians, and if $H$ is the hyperplane
section class we have $\Cal O_{F(n)}(H)=\Cal O_{F(n)}(1,1,\dots ,1)=
\Cal O_{F(n)}(\Omega_{s_1}+\dots +\Omega_{s_{n-1}})$. Also, the canonical
class is $K_{F(n)}=-2H$.(see \cite {M})

A map $f:\Bbb P^1\rightarrow F(V^{\ast })$ of multidegree
$\overline d=(d_1,d_2,\dots ,d_{n-1})$ with respect to $H$ is given by
specifying a flag of subbundles $S_1\subset S_2 \subset \dots
\subset S_{n-1}\subset V^{\ast } \otimes \Cal O_{\Bbb P^1}$ with
rank$(S_i)=i$, deg$(S_i)=-d_i$. Fixing the multidegree of the map amounts to
fixing the homology class $f_{\ast }[\Bbb P^1]=\sum d_iY_i$, where
$Y_i\in H_2(F(V^{\ast });\Bbb Z)$ is the Poincar\'e dual of $\Omega _{s_i}$.
Since $F(V^{\ast })$ is a homogeneous space, the moduli space of such maps
$Hom_{\overline d}(\Bbb P^1,F(V^{\ast }))$ is a smooth quasiprojective
variety of dimension $h^0(\Bbb P^1,f^{\ast }T_{F(V^{\ast })})=
\text {dim} F(V^{\ast })+f_{\ast }[\Bbb P^1]\cdot (-K_{F(V^{\ast })})=
\frac{n(n-1)}{2}+2\sum_{i=1}^{n-1}d_i$.

Our main tool is a certain compactification of $Hom_{\overline d}
(\Bbb P^1,F(V^{\ast }))$, generalizing Grothendieck's Quot scheme, which we
will introduce now.

Let $X$ be a smooth projective variety over an algebraically closed field $k$
and $E$ a vector bundle on $X$. For any scheme $S$ over $k$, let
$\pi :X\times S\rightarrow X$ be the projection. Consider the functor
$\Cal F(X,E):\{ Schemes\ over\ k\} \rightarrow \{ Sets\}$ given by
$$ \Cal F(X,E)(S)=\left \{ \matrix \text {equivalence classes of flagged
quotient sheaves} \\ \pi ^{\ast }E\twoheadrightarrow Q_{n-1}
\twoheadrightarrow \dots \twoheadrightarrow Q_1, \text{which are flat over}\
S \endmatrix \right \} ,$$
and for a morphism $S\overset \varphi \to \longrightarrow T$,
$$\Cal F(X,E)(\varphi)=\text { pull-back by}\ \varphi.$$
Here $ \pi ^{\ast }E\twoheadrightarrow Q_{n-1}\twoheadrightarrow
\dots \twoheadrightarrow Q_1$ is equivalent to $ \pi ^{\ast }E
\twoheadrightarrow Q'_{n-1}\twoheadrightarrow \dots \twoheadrightarrow Q'_1$
if there are isomorphisms $\theta _i:Q_i\rightarrow Q'_i$ such that all the
squares commute.

Let ${\overline P}=(P_1(m),\dots ,P_{n-1}(m))$ be numerical polynomials and
define the subfunctor $\Cal F_{\overline P}(X,E)$ by requiring that
$\chi (X_s,Q_{i_s})=P_i(m)$ for all $s\in S$. Extending the construction of
the Quot scheme we have the following
\proclaim{Theorem 1} For fixed ${\overline P}(m)$, $\Cal F_{\overline P}(X,E)$
 is represented by a projective scheme.
\endproclaim
 We will denote this scheme by $\Cal H\Cal Q_{\overline P}(X,E)$ and refer to
it as the {\it hyper-quot scheme} associated to $X$, $E$ and ${\overline P}$.
In general this scheme may be very complicated, but in the case of interest
for our purposes it is well-behaved. More precisely, let $P_i(m)=(m+1)i+d_i$
which is the Hilbert polynomial of a vector bundle of rank $i$ and degree
$d_i$ on $\Bbb P^1$. Then ${\overline P}$ is determined by
${\overline d}=(d_1,\dots ,d_{n-1})$
only. Denote the hyper-quot scheme associated to $\Bbb P^1$, $ V^{\ast }
 \otimes \Cal O_{\Bbb P^1}$ and ${\overline d}$ by $\Cal H\Cal Q_{\overline d}
(\Bbb P^1,F(V^{\ast }))$.
\proclaim{Theorem 2} $\Cal H\Cal Q_{\overline d}(\Bbb P^1,F(V^{\ast }))$ is a
smooth projective variety of dimension $\frac{n(n-1)}{2}+2\sum_{i=1}^{n-1}
d_i$, containing $Hom_{\overline d}(\Bbb P^1,F(V^{\ast }))$ as an open
subscheme. \endproclaim
 As a fine moduli space, $\Cal H\Cal Q_{\overline d}( \Bbb P^1,F(V^{\ast }))$
comes equipped with an universal sequence of sheaves $$ 0\hookrightarrow
\Cal S_1\hookrightarrow \Cal S_2\hookrightarrow \dots \hookrightarrow
\Cal S_{n-1}\hookrightarrow \Cal S_n=V^{\ast } \otimes \Cal O
\twoheadrightarrow \Cal T_{n-1}\twoheadrightarrow \dots \twoheadrightarrow
\Cal T_2\twoheadrightarrow \Cal T_1\twoheadrightarrow 0$$ on $\Bbb P^1\times
\Cal H\Cal Q_{\overline d}(\Bbb P^1,F(V^{\ast }))$, where $\Cal S_i=\text
{ker} \{ V^{\ast } \otimes \Cal O \twoheadrightarrow \Cal T_{n-i}\} $. For
$i=1,\dots ,n-1,\ \Cal S_i$ is a vector bundle of rank $i$ and relative degree
$-d_i$ (this follows from flatness), but the inclusions are injective as maps
of {\it sheaves} only!

 In fact, $Hom_{\overline d}(\Bbb P^1,F(V^{\ast }))$ is the largest subscheme
$U$ of $\Cal H\Cal Q_{\overline d}( \Bbb P^1,F(V^{\ast }))$ with the property
that on $\Bbb P^1\times U$ all the inclusions are injections of vector
bundles. The main technical result needed for our computations describes the
locus where these maps degenerate. Let $\overline e=(e_1,\dots ,e_{n-1})$
 with $e_i\leq max(i,d_i)$ and $\overline e_i=(0,\dots ,0,1,0,\dots ,0)$.
\proclaim{Theorem 3} There are rational maps $j_{\overline d-\overline e}:
\Bbb P^1\times \Cal H\Cal Q_{\overline d-\overline e}( \Bbb P^1,F(V^{\ast }))
--\rightarrow \Cal H\Cal Q_{\overline d}( \Bbb P^1,F(V^{\ast }))$ with the
following properties:

(i) $j_{\overline d-\overline e}$ is defined on $\Bbb P^1 \times
Hom_{\overline d-\overline e}(\Bbb P^1,F(V^{\ast }))$
and its restriction to it is an embedding.

(ii) If $s\in  \Cal H\Cal Q_{\overline d}( \Bbb P^1,F(V^{\ast }))$ is in the
image of $j_{\overline d-\overline e}$, then $\operatorname {rank}(\Cal S_i
\hookrightarrow \Cal S_{i+1})\leq i-e_i$ at (t,s), for some $t\in \Bbb P^1.$

 (iii) The images of the maps $j_{\overline d-\overline e}$ cover
$\Cal H\Cal Q_{\overline d}( \Bbb P^1,F(V^{\ast }))\setminus Hom_{\overline d}
(\Bbb P^1,F(V^{\ast }))$.
\endproclaim
{}From the previous theorem, one can see that the boundary of the hyper-quot
scheme compactifying the moduli space of maps $\Bbb P^1\rightarrow
F(V^{\ast })$ consists of $n-1$ divisors
$\Bbb D_1,\dots ,\Bbb D_{n-1}$, which are birational to $\Bbb P^1\times
\Cal H\Cal Q_{\overline d-\overline e_1}( \Bbb P^1,F(V^{\ast })),\dots ,
\Bbb P^1\times \Cal H\Cal Q_{\overline d-\overline e_{n-1}}( \Bbb P^1,
F(V^{\ast })). $

Following \cite {Be2}, we will define now the quantum multiplication
for cohomology classes on
$F(V^{\ast })$. There is an evaluation morphism $$ev: \Bbb P^1\times
Hom_{\overline d}(\Bbb P^1,F(V^{\ast })) \rightarrow F(V^{\ast }),$$ given by
$ ev(t,f)=f(t)$.

For $t\in \Bbb P^1$ , $w\in S_n$, define a subscheme of
$Hom_{\overline d}(\Bbb P^1,F(V^{\ast }))$ by $$\Omega_w(t)=
ev^{-1}(\Omega _w)\bigcap\{ \{t\} \times  Hom_{\overline d}
(\Bbb P^1,F(V^{\ast }))\}.$$ Set theoretically,$$\Omega_w(t)=\{ f\in
Hom_{\overline d}(\Bbb P^1,F(V^{\ast }))\ |\ f(t)\in \Omega_w\} .$$ Also, we
define $\overline{\Omega}_w(t)$ to be the following subscheme of
$\Cal H\Cal Q_{\overline d}(\Bbb P^1,F(V^{\ast }))$: consider the dual
sequence $$V \otimes \Cal O\rightarrow \Cal S_{n-1}^{\ast }\rightarrow \dots
\rightarrow \Cal S_1^{\ast }$$ on $\Bbb P^1\times \Cal H\Cal Q_{\overline d}
(\Bbb P^1,F(V^{\ast }))$ and the fixed flag $0=V_0\subset V_1\subset \dots
\subset V_{n-1}\subset V_n=V$. Let $$D_w^{p,q}=\text {locus where rank}
(V_p\otimes \Cal O \rightarrow \Cal S_q^{\ast })\leq r_w(q,p)$$ and $$D_w^{p,q}
(t)=D_w^{p,q}\bigcap \{ \{t\} \times \Cal H\Cal Q_{\overline d}(\Bbb P^1,
F(V^{\ast }))\}.$$ Then
$$\overline{\Omega}_w(t):=\bigcap_{p,q=1}^{n-1}D_w^{p,q}(t).$$

The following lemma and its corollaries will allow us to define
intersection numbers on $Hom_{\overline d}(\Bbb P^1,F(V^{\ast }))$. The proof
of $(ii)$ of the lemma depends heavily on the analysis in Theorem 3.
\proclaim{(Moving) Lemma} (i) For any $w_1,\dots ,w_N\in S_n$;$\ t_1,\dots
,t_n\in \Bbb P^1$ and general translates of
$\Omega _{w_i}\subset F(V^{\ast })$, the
intersection $\bigcap_{i=1}^N \Omega_{w_i}(t_i)$ has pure codimension
$\sum_{i=1}^N \ell (w_i)$ in $Hom_{\overline d}(\Bbb P^1,F(V^{\ast }))$ (or
is empty).

(ii) If $t_1,\dots ,t_N$ are $\underline {distinct}$, then for general
translates of the $\Omega_{w_i}$ the intersection $\bigcap_{i=1}^N
\overline {\Omega }_{w_i}(t_i)$ has pure codimension $\sum_{i=1}^N \ell (w_i)$
in $\Cal H\Cal Q_{\overline d}(\Bbb P^1,F(V^{\ast }))$ and is the (Zariski)
closure of $\bigcap_{i=1}^N \Omega_{w_i}(t_i)$ (or is empty).
\endproclaim
\proclaim{Corollary 1} The class of $\overline {\Omega}_w(t)$ in
$CH^{\ell (w)}(\Cal H\Cal Q_{\overline d}(\Bbb P^1,F(V^{\ast })))$ is
independent of $t\in \Bbb P^1$and the flag $V_{\bullet }\subset V$.
\endproclaim
Let $\Cal D=\frac {n(n-1)}{2}+2\sum_{i=1}^{n-1}d_i=\text {dim}
(\Cal H\Cal Q_{\overline d}(\Bbb P^1,F(V^{\ast })))$.
\proclaim{Corollary 2} If $\sum_{i=1}^N\ell (w_i)=\Cal D$, and
$t_1,\dots ,t_N$ are distinct, then the number of points in
$\bigcap_{i=1}^N \Omega_{w_i}(t_i)$ can be computed as the degree of
$(\bigcap_{i=1}^N \overline {\Omega }_{w_i}(t_i))$ in $CH^{\Cal D }
( \Cal H\Cal Q_{\overline d}(\Bbb P^1,F(V^{\ast })))$ (hence it is
independent of $t_i$ and the translates of $\Omega_{w_i}$).
\endproclaim
The corollaries imply that we have a well defined intersection number
$$\langle \Omega_{w_1},\dots ,\Omega_{w_N}\rangle _{\overline d}:
= \cases \text {number of points in} \bigcap_{i=1}^N \Omega_{w_i}(t_i),&\text
{if} \sum_{i=1}^N\ell (w_i)=\Cal D,\\  0,&\text { otherwise} . \endcases$$
This is the {\it Gromov-Witten} invariant associated to the classes
$\Omega_{w_1},\dots ,\Omega_{w_N}$.

{\bf Definition:} The quantum multiplication map is the linear map
$$m_q:\operatorname {Sym}(H^{\ast }(F(V^{\ast });\Bbb C)[q_1, \dots ,q_{n-1}])
\rightarrow H^{\ast }(F(V^{\ast });\Bbb C)[q_1, \dots ,q_{n-1}]$$ given by
$$m_q(\prod _{i=1}^N \Omega _{w_i}\overline q ^{\overline m _i})=
\overline q ^{\sum\overline m _i}\sum_{\overline d \in \Bbb N^{n-1}}
\overline q ^{\overline d}(\sum_{w\in S_n}\langle \Omega_w, \Omega_{w_1},
\dots ,\Omega_{w_N}\rangle _{\overline d}X_w).$$
Here $\overline q ^{\overline m}$ denotes as usual the monomial
$\prod _{i=1}^{n-1}q_i^{m_i}$.

 \cite {RT} proves that the pairing induced by $m_q$ on $\operatorname {Sym}^2
(H^{\ast }(F(V^{\ast });\Bbb C)[q_1,\dots ,q_{n-1}])$ determines a ring
structure on $H^{\ast }(F(V^{\ast });\Bbb C)[q_1, \dots ,q_{n-1}]$
( associativity of quantum multiplication; see also \cite {LT}.)
The pair ($H^{\ast }(F(V^{\ast });
\Bbb C)[q_1, \dots ,q_{n-1}],m_q)$ is the {\it quantum cohomology ring} of
$F(V^{\ast })$.
It is easy to prove (see e.g. \cite {Be2}), using the Moving Lemma, that
$m_q$ is the identity map on
$\operatorname {Sym}^1(H^{\ast }(F(V^{\ast });\Bbb C)[q_1, \dots ,q_{n-1}])$.
The restriction  $$m_q:\Bbb C[x_1,x_2, \dots ,x_n,q_1,q_2,\dots ,q_{n-1}]
\rightarrow H^{\ast }(F(V^{\ast });\Bbb C)[q_1, \dots ,q_{n-1}]$$ is
surjective (see \cite {ST} for an easy argument
by induction on degree) and by \cite {RT} it is a
ring homomorphism. Let $I$ be the kernel
of this map. Then the quantum cohomology ring of $F(n)$ is $$(H^{\ast }
(F(V^{\ast });\Bbb C)[q_1, \dots ,q_{n-1}],m_q)\cong \Bbb C[x_1,x_2, \dots ,
x_n,q_1,q_2,\dots ,q_{n-1}]/I.$$
We compute here the generators of $I$.

Recall that in $H^{\ast }(F(V^{\ast });\Bbb C)$ we have $R_k(n)=0$, where
$R_k(n)$ is the $k^{th}$ symmetric function in $x_1,\dots ,x_n,\ k=1,\dots ,n$.
 In the quantum ring however, $m_q(R_k(n))$ is a polynomial  $R'_k(n)(x_1,
\dots ,x_n,q_1,\dots ,q_{n-1})$ which doesn't vanish anymore (unless $k=1$).

{\bf Definition.} The {\it quantum deformation} of $R_k(n)$ is $ R^q_k(n):=
R_k(n)-R'_k(n).$
\proclaim{Theorem 4} (i) The quantum deformations of the relations $R_k(n)$
can be computed recursively with the formula $$R^q_k(n)=R^q_k(n-1)+x_n
\cdot R^q_{k-1}(n-1)+q_{n-1}\cdot R^q_{k-2}(n-2).$$
(Here $R^q_n(n-1)$ is set to be 0 and $R^q_0(n-2)$ is set to be 1).

(ii) The ideal $I$ is generated by $R^q_1(n),\dots ,R^q_n(n)$.
\endproclaim
{\bf Remarks: 1.} From Theorem 4 one can get a presentation for the quantum
cohomology ring of $F(n)$ for all $n>1$, once the ring is known for $F(1)$.
But $F(1)\cong \Bbb P^1$ and for projective spaces quantum cohomology is
well-known. In this case the ring is isomorphic to
$\Bbb C[x_1,x_2,q_1]/(x_1+x_2,\ x_1x_2+q_1)$.

{\bf 2.} In \cite {GK}, Givental and Kim gave the following (conjectural at
that time) compact description of the generators for the ideal $I$:

$R^q_k(n)$ is the coefficient of $\lambda ^{n-k}$ in the expansion of the
determinant: $$ \vmatrix x_1+\lambda & q_1 & 0 & 0 & \hdots & 0 & 0\\-1 &
 x_2+\lambda & q_2 & 0 & \hdots & 0 & 0\\ 0 & -1 & x_3+\lambda & q_3 & \hdots
 & 0 & 0\\ \vdots & \vdots & \vdots & \vdots & \vdots & \vdots & \vdots\\ 0 &
 0 & 0 & 0 & \hdots & x_{n-1}+\lambda & q_{n-1}\\ 0 & 0 & 0 & 0 & \hdots & -1
 & x_n \endvmatrix .$$
If one expands this determinant along the last column, the formula of Theorem
4 $(i)$ is obtained.

\heading Sketch of Proofs \endheading
For Theorem 1 and Theorem 2 the proofs, although lengthy, are straightforward
following the same line as for the corresponding results in the case of the
Quot functor.

For the proof of Theorem 3 we will change notation slightly and let
$\Cal S_i^{\overline d-\overline e}$ be the $i^{th}$ universal bundle on
$\Bbb P^1\times \Cal H\Cal Q_{\overline d-\overline e}( \Bbb P^1,
F(V^{\ast }))$.

 The maps $j_{\overline d-\overline e}$ can be constructed by downward
recursion. Therefore, we will outline here the construction for
$j_{\overline d-\overline e_i}, i=1, \dots ,n-1$. Let
$$\Bbb P^1\times \Bbb P^1\times \Cal H\Cal Q_{\overline d-\overline e_1}
( \Bbb P^1,F(V^{\ast })) \overset \pi_1\to \longrightarrow \Bbb P^1\times
\Cal H\Cal Q_{\overline d-\overline e_1}( \Bbb P^1,F(V^{\ast }))$$
be the projection. If we denote by $\Delta$ the diagonal in $ \Bbb P^1\times
\Bbb P^1$ , then $$\Delta \times \Cal H\Cal Q_{\overline d-\overline e_1}
( \Bbb P^1,F(V^{\ast })) \subset \Bbb P^1\times \Bbb P^1\times \Cal H\Cal Q_{
\overline d-\overline e_1}( \Bbb P^1,F(V^{\ast }))$$ is a divisor. Pull-back
by $\pi_1$ of the universal sequence on $\Bbb P^1\times
\Cal H\Cal Q_{\overline d-\overline e_1}( \Bbb P^1,F(V^{\ast }))$
gives the following sequence on $\Bbb P^1\times \Bbb P^1\times \Cal H\Cal Q_{
\overline d-\overline e_1}( \Bbb P^1,F(V^{\ast }))$:
$$0\hookrightarrow \pi_1^{\ast }\Cal  S_1^{\overline d-\overline e_1}
\hookrightarrow \pi_1^{\ast }\Cal  S_2^{\overline d-\overline e_1}
\hookrightarrow \dots \hookrightarrow \pi_1^{\ast }\Cal  S_{n-1}^{\overline
d-\overline e_1}\hookrightarrow \pi_1^{\ast }\Cal  S_n^{\overline
d-\overline e_1} =V^{\ast } \otimes \Cal O_{\Bbb P^1\times\Bbb P^1\times
\Cal H\Cal Q_{\overline d-\overline e_1}}$$
Let $\Cal {\tilde S}_1^{\overline d-\overline e_1}=\bigl(\pi_1^{\ast }
\Cal  S_1^{\overline d-\overline e_1}\bigr)\otimes \Cal O \bigl(-\Delta
\times \Cal H\Cal Q_{\overline d-\overline e_1}( \Bbb P^1,F(V^{\ast }))\bigr)$.
 Then $\Cal {\tilde S}_1^{\overline d-\overline e_1}$ is a vector bundle of
rank 1 and relative degree $-d_1$ on $\Bbb P^1\times \Bbb P^1\times
\Cal H\Cal Q_{\overline d-\overline e_1}( \Bbb P^1,F(V^{\ast }))$, which is a
subsheaf of $\pi_1^{\ast }\Cal  S_2^{\overline d-\overline e_1}$. Since
$\Cal H\Cal Q_{\overline d}( \Bbb P^1,F(V^{\ast }))$ is a fine moduli space,
 we get a morphism
$$j_{\overline d-\overline e_1}:\Bbb P^1\times
\Cal H\Cal Q_{\overline d-\overline e_1}( \Bbb P^1,F(V^{\ast }))\longrightarrow
\Cal H\Cal Q_{\overline d}( \Bbb P^1,F(V^{\ast }))$$ such that
$$(id,j_{\overline d-\overline e_1})^{
\ast } \Cal S_i^{\overline d}=\cases \Cal {\tilde S}_1^{\overline d-\overline
e_1},&\text{for $i=1$}\\ \pi_1^{\ast }\Cal  S_i^{\overline d-\overline e_1},&
\text{for $i\neq 1$}.\endcases$$
It is easy to see that if rank$\{\Cal S_1^{\overline d}\hookrightarrow
\Cal S_2^{\overline d}\}=0$ at $(t,x)\in \Bbb P^1\times \Cal H\Cal Q_{
\overline d}$, then $x$ is in the image of $j_{\overline d-\overline e_1}$
and that the restriction of $j_{\overline d-\overline e_1}$ to $\Bbb P^1\times
Hom_{\overline d-\overline e_1}( \Bbb P^1,F(V^{\ast }))$ is an embedding.

For $2\leq i\leq n-1$, let
$$\Bbb P^1\times \Bbb P^1\times \Cal H\Cal Q_{\overline d-\overline e_i}
( \Bbb P^1,F(V^{\ast })) \overset \pi_i\to \longrightarrow \Bbb P^1\times
\Cal H\Cal Q_{\overline d-\overline e_i}( \Bbb P^1,F(V^{\ast }))$$
be the projection and consider as before the sequence
$$0\hookrightarrow \pi_i^{\ast }\Cal  S_1^{\overline d-\overline e_i}\dots
\hookrightarrow \pi_i^{\ast }\Cal  S_{i-1}^{\overline d-\overline e_i}
\hookrightarrow \pi_i^{\ast }\Cal  S_i^{\overline d-\overline e_i}\dots
\hookrightarrow \pi_i^{\ast }\Cal  S_{n-1}^{\overline d-\overline e_i}
\hookrightarrow V^{\ast } \otimes \Cal O_{\Bbb P^1\times\Bbb P^1\times
\Cal H\Cal Q_{\overline d-\overline e_i}}.$$
Let $\Cal L_i^{\overline d-\overline e_i}$ be the quotient
$\pi_i^{\ast }\Cal  S_i^{\overline d-\overline e_i}/
\pi_i^{\ast }\Cal S_{i-1}^{\overline d-\overline e_i}$.
The map $\Cal L_i^{\overline d-\overline e_i}(-\Delta
\times \Cal H\Cal Q_{\overline d-\overline e_i})
\longrightarrow \Cal L_i^{\overline d-\overline e_i}$ induces
$$\operatorname {Ext}^1\bigl(\Cal L_i^{\overline d-\overline e_i},
\pi_i^{\ast }\Cal  S_{i-1}^{\overline d-\overline e_i}\bigr) \longrightarrow
\operatorname {Ext}^1\bigl(\Cal L_i^{\overline d-\overline e_i}(-\Delta
\times \Cal H\Cal Q_{\overline d-\overline e_i}),
\pi_i^{\ast }\Cal  S_{i-1}^{\overline d-\overline e_i}\bigr).$$
Thus we have a diagram of sheaves on $\Bbb P^1\times\Bbb P^1\times
\Cal H\Cal Q_{\overline d-\overline e_i}( \Bbb P^1,F(V^{\ast }))$:
$$\CD
0 @>>> \pi_i^{\ast }\Cal  S_{i-1}^{\overline d-\overline e_i}
 @>>>  \pi_i^{\ast }\Cal  S_i^{\overline d-\overline e_i}
 @>>> \Cal L_i^{\overline d-\overline e_i} @>>> 0 \\
@.     @|   @AAA   @AAA   @. \\
0 @>>> \pi_i^{\ast }\Cal  S_{i-1}^{\overline d-\overline e_i}
 @>>>  \Cal {\tilde S}_i^{\overline d-\overline e_i}
 @>>> \Cal L_i^{\overline d-\overline e_i}(-\Delta
\times \Cal H\Cal Q_{\overline d-\overline e_i})
 @>>> 0 \endCD$$
The restriction of $\Cal {\tilde S}_i^{\overline d-\overline e_i}$ to the open
set $\pi_i^{-1}(\Cal U)=\Bbb P^1 \times \Cal U$, where $$\Cal U=\{y\in
\Bbb P^1\times \Cal H\Cal Q_{\overline d-\overline e_i} \mid \text {rank}_y\{
\Cal S_{i-1}^{\overline d-\overline e_i}\hookrightarrow \Cal  S_i^{\overline
d-\overline e_i}\}=i-1\},$$ is a vector bundle of rank $i$ and relative degree
$-d_i$. Also, the restriction of
$$\Cal {\tilde S}_i^{\overline d-\overline e_i}\longrightarrow
\pi_i^{\ast }\Cal  S_i^{\overline d-\overline e_i}$$
to $\Bbb P^1 \times \Cal U$ is injective as a map of sheaves. Note that
$\Bbb P^1\times Hom_{\overline d-\overline e_i}( \Bbb P^1,F(V^{\ast }))
\subset \Cal U$.

It follows that there is a rational map
 $$j_{\overline d-\overline e_i}:\Bbb P^1\times
\Cal H\Cal Q_{\overline d-\overline e_i}( \Bbb P^1,F(V^{\ast }))-\ -\
\rightarrow \Cal H\Cal Q_{\overline d}( \Bbb P^1,F(V^{\ast }))$$
which is defined on $\Cal U$ and has the following properties:

1. On $\Bbb P^1 \times \Cal U$,
$$(id,j_{\overline d-\overline e_i})^{
\ast } \Cal S_k^{\overline d}=\cases \Cal {\tilde S}_i^{\overline d-\overline
e_i},&\text{for $k=i$}\\ \pi_i^{\ast }\Cal  S_k^{\overline d-\overline e_i},&
\text{for $k\neq i$}.\endcases$$

2. If $\Cal S_i^{\overline d}\hookrightarrow \Cal  S_{i+1}^{\overline
d}$ degenerates at $(t,x)\in \Bbb P^1 \times \Cal H\Cal Q_{\overline d}$ and
$\Cal S_{i-1}^{\overline d}\hookrightarrow \Cal  S_i^{\overline d}$ has maximal
rank at $(t,x)$, then $x$ is in the image of $j_{\overline d-\overline e_i}$.

3. The restriction of $j_{\overline d-\overline e_i}$ to $\Bbb P^1\times
Hom_{\overline d-\overline e_i}( \Bbb P^1,F(V^{\ast }))$ is an embedding.

It is clear that $$\Cal H\Cal Q_{\overline d}( \Bbb P^1,F(V^{\ast }))
\setminus Hom_{\overline d}(\Bbb P^1,F(V^{\ast }))=\bigcup_{i=1}^{n-1}
\operatorname {Im}(j_{\overline d-\overline e_i}).$$

{\bf Remark 3.} We can describe now the preimages of the degeneracy loci
$\overline \Omega_w(t)$ by the maps $j_{\overline d-\overline e_i}$. Note
first that the map $\Cal {\tilde S}_i^{\overline d-\overline e_i}
\longrightarrow \pi_i^{\ast }\Cal  S_i^{\overline d-\overline e_i}$
is an isomorphism outside $\Delta \times \Cal H\Cal Q_{\overline d-
\overline e_i}$ and that when restricted to $\Delta \times \Cal H\Cal Q_
{\overline d-\overline e_i}$ it factors trough $\pi_i^{\ast }\Cal  S_{i-1}^{
\overline d-\overline e_i}$. Consequently, $j_{\overline d-\overline e_i}^{
-1}(\overline \Omega_w(t))$ splits as the union of two subschemes of
$\Bbb P^1 \times \Cal H\Cal Q_{\overline d-\overline e_i}$:
$$j_{\overline d-\overline e_i}^{-1}(\overline \Omega_w(t))=\Bbb P^1\times
\overline \Omega_w(t)\ \bigcup \ {\tilde \Omega}_w(t),$$ where
${\tilde \Omega}_w(t)$ is the locus $\bigcap_{q\neq i} D_w^{p,q}(t)$ in
$\{ t\} \times \Cal H\Cal Q_{\overline d-\overline e_i}$.

Using the previous remark, the same proof as in \cite {Be2}, Lemma 2.2  and
2.2 A, gives the Moving Lemma in our case too.

We will give now in more detail the proof of Theorem 4.

The following result, due to Fulton (Theorem 8.2.
in \cite {F1}), will be crucial. It generalizes the Giambelli formula
mentioned in the introduction. We give here only the version needed for our
purposes and refer to the original paper for the full statement.

We need to establish some notations first. Let $X$ be a
scheme of finite type over a field $k$, $V$ an $n$-dimensional vector space,
$V_{\bullet }\subset V$ a fixed complete flag of subspaces and
$w\in S_{n+1}$ a permutation. Consider a flag of vector bundles on $X$:
$$B_{n-1}\twoheadrightarrow B_{n-2}\twoheadrightarrow \dots \twoheadrightarrow
B_1,$$ with rank$(B_i)=i$. To each map $h:V\otimes \Cal O_X \longrightarrow
B_{n-1}$ one can associate degeneracy loci $\Omega_w(h)$ defined by the
conditions
$$\text {rank}\{ V_p\otimes \Cal O_X \rightarrow B_q\} \leq r_w(q,p)$$
for $1\leq p\leq n$ and $1\leq q \leq n-1$. Let $x_1=c_1(B_1),\ x_i=c_1(\text
{ker}\{ B_i\rightarrow B_{i-1}\})$ for $i=2,\dots ,n-1$, $x_n=0$.
 \proclaim {Theorem (Fulton)} If X is a purely d-dimensional scheme, there
is a class $\hat {\Omega}_w(h)$ in $CH_{d-\ell (w)}(\Omega_w(h))$,
satisfying the following.

(i) The image of $\hat { \Omega}_w(h)$ in $CH_{d-\ell (w)}(X)$ is
$\sigma_w(x)\cap [X]$.

(ii) Each irreducible component of $\Omega_w(h)$ has codimension at most
$\ell (w)$ in X. If $\text {codim}(\Omega_w(h),X)=\ell (w)$, then
$\hat { \Omega}_w(h)$ is a positive cycle whose support is $\Omega_w(h)$.

(iii) If $\text {codim}(\Omega_w(h),X)=\ell (w)$ and X is Cohen-Macaulay, then
$\Omega_w(h)$ is Cohen-Macaulay and $$\hat { \Omega}_w(h)=[\Omega_w(h)].$$
  \endproclaim
{\bf Remark 4.} It will be important later that the statement is formulated
for $w\in S_{n+1}$ rather than $S_n$. Of course, the theorem applies to the
degeneracy loci we employed so far, since both the loci and the Schubert
polynomials are invariant under the natural embedding $S_n\subset S_m,\
\text {for}\ m>n$.

By the previous theorem, if the maps
$\Cal S_i^{\ast }\rightarrow \Cal S_{i-1}^{\ast }$
were {\it surjective} for $i=2,\dots ,n-1$, one could express the degeneracy
locus $\overline {\Omega}_w(t)$ in $CH^{\ast }(\Bbb P^1\times
\Cal H\Cal Q_{\overline d}( \Bbb P^1,F(V^{\ast })))$ as the Schubert
polynomial $\Cal \sigma _w$ evaluated in the Chern classes $c_1(\Cal S_1^{
\ast }), c_1(\text {ker}\{ \Cal S_2^{\ast }\rightarrow\Cal S_1^{\ast }\} ), c_1
(\text {ker}\{ \Cal S_{n-1}^{\ast }\rightarrow\Cal S_{n-2}^{\ast }\} )$. In
our case, this is true only if we restrict to the open set $U$ of $\Bbb P^1
\times \Cal H\Cal Q_{\overline d}( \Bbb P^1,F(V^{\ast }))$ where $\Cal S_i$
is a {\it subbundle} of $\Cal S_{i+1}$ for all $i=1,\dots ,n-2$.

 Hence, if we let
$$\align x_i(t)&=c_1( \Cal S_i^{\ast })\bigcap \{ \{t\} \times
\Cal H\Cal Q_{\overline d}(\Bbb P^1,F(V^{\ast }))\}-c_1( \Cal S_{i-1}^{\ast })
\bigcap \{ \{t\} \times \Cal H\Cal Q_{\overline d}(\Bbb P^1,F(V^{\ast }))\}\\
&=\overline \Omega _{s_i}(t)-\overline \Omega _{s_{i-1}}(t),\endalign$$ we
have in $CH_{\ast }( \Cal H\Cal Q_{\overline d}( \Bbb P^1,F(V^{\ast })))$ the
following identity:
$$\overline \Omega_w(t)-\Cal \sigma _w(x_1(t),\dots ,x_{n-1}(t))=
j_{\ast }(G_w(t)),\tag 1 $$ where $G_w(t)\in CH_{\ast }(\bigcup _{i=1}^{n-2}
\Bbb D_i)$ and $j:\bigcup _{i=1}^{n-2}\Bbb D_i\rightarrow
\Cal H\Cal Q_{\overline d}( \Bbb P^1,F(V^{\ast }))$ is the inclusion.

Let $w_1,\dots ,w_N\in S_n$ such that $\ell(w)+\sum _{i=1}^N
\ell(w_i)=\Cal D$ and $t,t_1,\dots ,t_N\in \Bbb P^1$ distinct points. Then we
get$$\overline \Omega_w(t)\cdot \prod _{i=1}^N
\overline {\Omega }_{w_i}(t_i)-\Cal \sigma _w(x_1(t),\dots ,
x_{n-1}(t))\cdot \prod _{i=1}^N\overline {\Omega }_{w_i}(t_i)= \tag 2$$
$$j_{\ast }(G_w(t))\cdot \prod _{i=1}^N
\overline {\Omega }_{w_i}(t_i)$$ in $CH^{\ast }( \Cal H\Cal Q_{\overline d}
( \Bbb P^1,F(V^{\ast })))$. From the Moving Lemma and its corollaries, we see
that $j_{\ast }(G_w(t))\cdot \prod _{i=1}^N\overline {\Omega }_{w_i}(t_i)$
 is (the negative of) the (signed) number of points in
$\Cal \sigma _w(x_1(t),\dots ,x_{n-1}(t))\cdot \prod _{i=1}^N\overline
{\Omega }_{w_i}(t_i)$ supported on $\bigcup _{i=1}^{n-2}\Bbb D_i$.

Note that for distinct points $u_1, \dots ,u_{n-1},t_1,\dots ,t_N\in \Bbb P^1$
the intersection $$\Cal \sigma _w(\overline {\Omega}_1(u_1),\dots ,
\overline {\Omega}_{n-1}(u_{n-1})-
\overline {\Omega}_{n-2}(u_{n-2}))\bigcap \bigcap _{i=1}^N\overline
{\Omega }_{w_i}(t_i)$$ avoids the boundary. Allowing the points $u_1, \dots
,u_{n-1}$ to come together ``moves'' part of the intersection in the boundary.

For $m\leq n,\ k\leq m-1$, let $\alpha _k(m)\in S_n$ be the permutation
$$\pmatrix 1&2&\dots & m-k-1&m-k&\dots &m-1&m&\dots &n-1&n  \\1&2&\dots &m-k-1
&m-k+1&\dots &m&m-k&\dots &n-1&n \endpmatrix .$$ Its Schubert polynomial is
the $k-1^{st}$ elementary symmetric function in variables $x_1,\dots
,x_{m-1}$: $\ \Cal \sigma_{\alpha _k(m)}(x)=R_{k-1}(m-1)$.

For the proof of the formula in Theorem 4 $(i)$ we will need to compute the
number given by the LHS of (2) for $w=\alpha _k(m)$.

Note first that $\Cal \sigma_{\alpha _k(m)}(x)=R_{k-1}(m-1)$ is the sum of all
degree $k-1$ ``square-free'' monomials in $x_1,\dots ,x_{m-1}$. Recall that
$x_p(t)=\overline \Omega _{s_p}(t)-\overline \Omega _{s_{p-1}}(t)$.
We are therefore led to compute the part supported on $\bigcup _{i=1}^{n-2}
\Bbb D_i$ for intersections of the type
$$\overline {\Omega }_{s_{i_1}}(t)\bigcap \dots \bigcap
\overline {\Omega }_{s_{i_{k-1}}}(t)\bigcap \bigl( \bigcap_{i=1}^N
\overline {\Omega }_{w_i}(t_i)\bigr),$$
with $i_j\in \{1, \dots ,m-1\}$. The $i_j$'s are not necessarily distinct, but
each index in $\{1, \dots ,m-2\}$ may occur at most twice, while $m-1$ may
occur at most once.

By Remark 3,
$$j_{\overline d-\overline e_i}^{-1}(\overline \Omega_{s_p}(t))=\cases
\Bbb P^1\times \overline \Omega_{s_p}(t)\ \bigcup \{ t\}\times \Cal H
\Cal Q_{\overline d-\overline e_p},&\text {for $i=p$}\\ \Bbb P^1\times
\overline \Omega_{s_p}(t)\bigcup \{ t\}\times \overline \Omega_{s_p}(t),&
\text{for $p\neq i$} \endcases$$ in $\Bbb P^1\times
\Cal H\Cal Q_{\overline d-\overline e_i},\ \text{for} \ i=1, \dots ,n-2$.

Since the points $t,t_1, \dots ,t_N \in \Bbb P^1$ are distinct, the only time
we will get nonempty intersections is when at least one index occurs twice.
Say that $\overline \Omega_{s_p}(t)$ occurs twice, for some $1\leq p\leq m-2$
(this is equivalent to saying that the product $x_px_{p+1}$ is in the
corresponding square-free monomial M). Then the intersection will be
$$\overline {\Omega }_{s_{l_1}}(t)\bigcap \dots \bigcap
\overline {\Omega }_{s_{l_{k-3}}}(t)\bigcap \bigl( \bigcap_{i=1}^N
\overline {\Omega }_{w_i}(t_i)\bigr) \subset \{ t\} \times \Cal H\Cal Q_{
\overline d-\overline e_p},$$
where the collection of indices $l_1,\dots ,l_{k-3}$ is obtained from $i_1,
\dots ,i_{k-1}$ by removing $ p$ (twice). This
corresponds to the monomial $M/x_px_{p+1}$. If $l_1,\dots ,l_{k-3}$ are not
distinct we can repeat the procedure, working now with
$\Cal H\Cal Q_{\overline d-\overline e_p}$. This says the following: in the
quantum ring of $F(n)$, for all $m\leq n$, we have
$$m_q(\Omega _{\alpha _k(m)},\bullet )=m_q(\Cal \sigma_{\alpha _k(m)}(x),
\bullet)+m_q(\Cal \sigma_{\alpha _k(m)}^\prime (x,q),\bullet),\tag 3 $$
where $\Cal \sigma_{\alpha _k(m)}^\prime (x,q)$ is the sum of all
square-free monomials of weighted
degree $k-1$ in $x_1,\dots ,x_{m-1},q_1,\dots ,q_{m-2}$
 (here deg$(x_i)=1$ and  deg$(q_i)=2$), with the
additional condition that if $q_i$ occurs in such a monomial, then none of
$x_i,x_{i+1},q_{i+1}$ can occur.

It is clear that $m_q(\Cal \sigma_{\alpha _k(m)}(x)+\Cal \sigma_{\alpha _k
(m)}^\prime (x,q))$ satisfies the recursion relation in Theorem 4$(i)$, i.e.
$$\split m_q(\Cal \sigma_{\alpha _k(m)}(x)+\Cal \sigma_{\alpha _k
(m)}^\prime (x,q))&=m_q(\Cal \sigma_{\alpha _k(m-1)}(x)+\Cal \sigma_{\alpha
_k(m-1)}^\prime (x,q)) \\ &+m_q(x_{m-1}(\Cal \sigma_{\alpha _{k-1}(m-1)}(x)+
\Cal \sigma_{\alpha _{k-1}(m-1)}^\prime (x,q))) \\ &+m_q(q_{m-2}(\Cal \sigma_{
\alpha _{k-2}(m-2)}(x)+\Cal \sigma_{\alpha _{k-2}(m-2)}^\prime (x,q))).
\endsplit \tag 4$$

{\bf Remark 5.} Recently, W. Fulton informed us that he computed the
coefficients of the characteristic polynomial obtained by expanding the
determinant in \cite {GK}, i.e. the generators of the ideal $I$. His formula,
which is easy to prove by induction, is (with our notations):
$$R_k(n)=\sigma_{\alpha_{k+1}(n+1)}(x)+
\sigma_{\alpha_{k+1}(n+1)}^\prime (x,q).$$
This suggested us to make the explicit calculation in (3), rather than showing
the recursion (4), which is the only thing we need to complete our proof.

 Monk's formula (see \cite {M}) gives an expression for the intersection of
Schubert subvarieties on $F(n)$: $$ \Omega _{s_p}\cdot \Omega _w=
\sum _{t_{ij}}\Omega_{w\cdot t_{ij}},$$ where the sum is over all
transpositions $t_{ij}$
of integers $i\leq p<j$ such that $\ell (w\cdot t_{ij})=\ell (w)+1$. As noted
in \cite {LS2}, this is not an identity among the corresponding Schubert
polynomials, unless one embedds $S_n$ in $S_{n+1}$ (as usual, by setting
$w(n+1)=n+1$). Using Monk's formula to multiply $\Cal \sigma_{s_{n-1}}(x)=
x_1+x_2+\dots +x_{n-1}$ and $\Cal \sigma _{\alpha _{k-1}(n)}(x)$ one sees
easily that we get $$ (x_1+x_2+\dots +x_{n-1})\cdot \Cal
\sigma _{\alpha _{k-1}(n)}(x)=\Cal \sigma_{\alpha _k(n)}(x)+\Cal
\sigma _{\beta _k}(x),\tag 5 $$ with $\beta _k=\alpha _{k-1}(n)
\cdot t_{n-1,n+1}\in S_{n+1}.$

 We can also define the degeneracy locus $\overline \Omega _{\beta _k}(t)\in
\Cal H\Cal Q_{\overline d}( \Bbb P^1,F(V^{\ast }))$ in a similar manner by
$$\overline \Omega _{\beta _k}(t)=\bigcap_{p=1}^n\bigcap_{q=1}^{n-1}
 D_w^{p,q}(t).$$
It is easy to see that $\overline \Omega _{\beta _k}(t)$ is given by
the conditions $\text {rank}(V^n\otimes \Cal O\rightarrow\Cal S_i^{\ast })
\leq i-1,\ \text {for}\ n-k+1\leq i\leq n-1$. In particular, it is supported
on $\Bbb D_{n-1}$.

 It is immediate from the definition that $x_1+\dots +x_n=0$ in the quantum
ring. Using (5), we obtain for $k\geq 2$ that
$$\split m_q(R_k(n))&=m_q(R_k(n-1)+x_nR_{k-1}(n-1))\\
&= m_q(R_k(n-1)-(x_1+x_2+\dots +x_{n-1})\cdot R_{k-1}(n-1))\\ &=
- \sum _{\overline d\in \Bbb N_{n-1}}\overline q^{\overline d}
\sum _{w\in S_n}\langle \Omega_w,\Cal \sigma_{\beta _k}(x_1,x_2,
\dots ,x_{n-1})\rangle _{\overline d}X_w.\endsplit \tag 6 $$
Now, by our definition and Fulton's Theorem 8.2 again,
$$\split  -\langle \Omega_w,\Cal \sigma_{\beta _k}(x_1,x_2,\dots ,x_{n-1})
\rangle _{\overline d}&=- \overline \Omega _w(u)\cdot \Cal \sigma_{\beta _k}
(x_1(t),\dots x_{n-1}(t))\\ &= \overline \Omega _w(u)\cdot j_{\ast }
(G_{\beta _k}(t))-\overline \Omega _w(u)\cdot \overline \Omega _{\beta _k}(t),
\endsplit \tag 7 $$ for distinct points $t,u\in \Bbb P^1$. (The same proof as
in the Moving Lemma shows that the codimension of $\overline
\Omega _{\beta _k}(t)$ in $\Cal H\Cal Q_{\overline d}( \Bbb P^1,F(V^{\ast }))$
is equal to $\ell (\beta _k))$.

 The description of $j_{\overline d-\overline e_{n-1}}^{-1}
(\overline \Omega _{\beta _k}(t))$ in Remark 3 gives, for
any $w\in S_n$,
$$\split \overline \Omega _w(u)\cdot \overline \Omega _{\beta _k}(t)&
= \overline \Omega _w(u)\cdot \overline \Omega _{\alpha _{k-2}(n-1)}(t) \\ &=
\overline \Omega _w(u)\cdot \Cal \sigma_{\alpha _{k-2}(n-1)}(x(t))+\overline
\Omega _w(u)\cdot j_{\ast }(G_{\alpha _{k-2}(n-1)}(t)),\endsplit \tag 8$$
where the first product is computed in $CH^{\ast }
( \Cal H\Cal Q_{\overline d}( \Bbb P^1,F(V^{\ast })))$ and the others are in
$CH^{\ast }( \Cal H\Cal Q_{\overline d-\overline e_{n-1}}( \Bbb P^1,
F(V^{\ast }))).$

Indeed, $\overline \Omega _{\beta _k}(t)=\bigcap_{q\neq n-1}
D_{\beta _k}^{p,q}$ and $r_{\alpha _{k-2}(n-1)}(q,p)=r_{\beta _k}(q,p)$
for all $1\leq p\leq n\ ,\ 1\leq q\leq n-2$, while $r_{\alpha _{k-2}(n-1)}
(n-1,p)=p$ for all $1\leq p\leq n$.

 Note that the Schubert polynomial of $\alpha _{k-1}(n-1)$ is
the $(k-2)^{nd}$ elementary symmetric function in $x_1,\dots ,x_{n-2}$
(and 1 if $k=2$).

Using (5) and Remark 3 again to compute the intersections supported on the
boundary, we get
$$\split \overline \Omega _w(u)\cdot j_{\ast }(G_{\beta _k}(t)) &=
-\overline \Omega _w(u)\cdot j_{\ast }(G_{\alpha _k(n)}(t))\\ &+
\overline \Omega _w(u)\cdot (x_1(t)+\dots +x_{n-1}(t))\cdot  j_{\ast }
(G_{\alpha _{k-1}(n)}(t))\endsplit . \tag 9$$

This is easy once we observe that $x_1(t)+\dots +x_{n-1}(t)=\overline
\Omega_{s_{n-1}}(t)$ and $$j_{\overline d-\overline e_i}^{-1}(\overline
\Omega_{s_{n-1}}(t))=\Bbb P^1 \times \overline \Omega_{s_{n-1}}(t)\
\bigcup\ \{ t\} \times \overline \Omega_{s_{n-1}}(t).$$
Combining (3), (4), (6), (7), (8) and (9), the formula in Theorem 4 $(i)$
follows by induction on $n$.

The proof of $(ii)$ is standard (see \cite {Bat}, \cite {Be2}, \cite {MP},
\cite {ST}).

{\bf Note added:} After this work was completed, B. Kim informed us that he
obtained independently Theorems 1 and 2.

\NoBlackBoxes

\enddocument